\documentclass[12pt]{article}
\title{\bf On the information-theoretic approach to G\"odel's incompleteness 
theorem}
\author{Germano D'Abramo\\
{\small Istituto di Astrofisica Spaziale e Fisica Cosmica,}\\
{\small Area di Ricerca CNR Tor Vergata, Roma, Italy}\\
{\small E--mail: {\tt dabramo@rm.iasf.cnr.it}}}

\date{}

\begin{document}

\maketitle

\begin{abstract} 

In this paper we briefly review and analyze three published proofs of
Chaitin's theorem, the celebrated information-theoretic version of
G\"odel's incompleteness theorem. Then, we discuss our main perplexity
concerning a key step common to all these demonstrations.

\end{abstract}

\section{Introduction}

The notion of algorithmic (or program-size) complexity, suggested for the
first time in the 1960s, identifies the complexity of a binary string with
the size of the shortest computer program that computes that string as
output. In the 1970s, mathematician G.~J.~Chaitin made use of this notion
to reformulate in terms of information-theoretic arguments the first
theorem of G\"odel on the incompleteness of formal axiomatic systems.
Chaitin's main result states that a formal axiomatic system coded in a
computer program of $n$ bits is unable to prove that a binary sequence
exists with a complexity greater than $n$ bits, whereas it can be easily
proved that there are actually infinitely many strings with this property.
The proof of Chaitin's result is based on the logical paradox which would
arise if a program of $n$ bits were able to prove the existence, and
therefore to generate a specific string of a complexity greater than $n$
bits. In this paper we briefly analyze three published proofs of Chaitin's
theorem and argue that they rely on a questionable, {\em a priori}
paradoxical assumption which seems to invalidate from the beginning the
demonstration {\em ab absurdo}: as a matter of fact, the logical paradox
reached at the end of the proofs seems to have been tacitly inserted at
the beginning of the proofs and therefore only a tautological, useless
{\em absurdum} is obtained.

\section{A questionable aspect of Chaitin's approach}

In order to explain what represents, according to the author, a dubious
point in Chaitin's approach to the incompleteness of formal systems, we
give here a more detailed description of his result.

In Chaitin~\cite{sa} and Gardner~\cite{ga} a formal axiomatic system is
described as a program for an idealized computer, a routine of $c$ bits
which makes use of $n$ bits of axioms and rules of inference to generate
systematically and check all possible proofs in order of their size: first
all proofs of one bit, then all proofs of two bits, and so forth (this
kind of computation is known as the {\it British Museum algorithm}). After
having generated a proof, according to Chaitin~\cite{sa} and
Gardner~\cite{ga}, the routine tests whether it is the first one proving
that a specific binary sequence is of a complexity greater than the number
of bits in the program, $n+c$. When the routine finds such a proof, it
prints the specific binary sequence and then halts. Now, the logical
paradox arises quite clearly: a program of $n+c$ bits supposedly
calculates a number that no program its size should be able to calculate.

This is the way in which Chaitin~\cite{sa} and Gardner~\cite{ga} show why
a formal axiomatic system of $n+c$ bits is unable to prove the existence
of any specific number with algorithmic complexity greater than $n+c$
bits. Actually, there are infinitely many number of complexity greater
than $n+c$ bits. As a matter of fact, with a number of bits less than or
equal to $n+c$ it is possible to create at the most $2^{n+c+1} - 2$
different binary strings\footnote{That is, we add up the number of all
possible strings of one bit ($2$), that of all possible strings of two
bits ($2^2$), and so on up to $n+c$ bits: $$\sum_{i=1}^{n+c} 2^i = 2
+2^2+2^3+...+2^{n+c}=2^{n+c+1}-2.$$}, each one of them is a program that
might generate a number. Thus, at the most $2^{n+c+1} - 2$ numbers have
algorithmic complexity less than or equal to $ n+c$ bits, while there are
infinitely many integer numbers.

Chaitin~\cite{ieee, co} shows that there is a program of length nearly
equal to $\log_2 N + k$ bits to calculate a number which supposedly can
not be calculated by a program shorter than $N$ bits.  Here $k$ is the
dimension of the proof-checking routine, essentially the formal axiomatic
system, and $\log_2 N$ is nearly equal to the dimension of the binary
expansion of decimal number $N$. Since $\log_2 N + k$ becomes much smaller
than $N$ for sufficiently large $N$, we have again that such formal
axiomatic system is unable to prove the existence of any number of a
complexity greater than $N$ bits.

Here we quote a passage from Chaitin~\cite{co} in which the working scheme
of such program is described:

\begin{quote}
{\em 

You start running through all possible proofs in the formal axiomatic
system in size order.  You apply the proof-checking algorithm to each
proof.  And after filtering out all the invalid proofs, you search for the
first proof that a particular positive integer requires at least an N-bit
program.

The algorithm that I've just described is very slow but it is very simple.  
It's basically just the proof-checking algorithm, which is $k$ bits long,
and the number $N$, which is $\log_2 N$ bits long.  So the total number of
bits is just

$$\log_2 N + k$$ 

as was claimed. That concludes the proof of my incompleteness result that
I wanted to discuss with G\"odel.}
\end{quote}

Let us now come to what appears to be a rather unsatisfactory aspect of
the approaches mentioned above. For the sake of simplicity, we concentrate
on the last proof, but our point holds also for proofs~\cite{sa},
\cite{ga} and \cite{ieee} with small changes.

It is not clear how the sole decimal number $N$ is enough to the formal
axiomatic system to identify univocally and unequivocally the desired
proof.  Keeping in mind the argument in Chaitin~\cite{co}, let us divide
the main algorithm of $\log_2 N + k$ bits into two programs of $k_1$ and
$\log_2 N + k_2$ bits respectively, such that $k_1+k_2\simeq k$ (meaning
``not too much greater than $k$''), and such that the concatenation of the
first and the second program works as the main algorithm. The first one is
the {\em pure proof-checking algorithm}, namely that which generates all
strings in increasing order of size and checks their formal correctness
within the formal axiomatic system. The second one accomplishes the {\em
pure searching task}; it mechanically runs through all valid proofs
provided by the first algorithm and tries to find a particular binary
string. In its searching task it uses as a guide only what is coded in its
own $\log_2 N + k_2$ bits. Like every pure searching algorithm, this
one is able to find univocally and unequivocally only what is exactly
expressible or compressible in at most $\log_2 N + k_2$ bits, not more.  
A successful search occurs in the following two cases:

\begin{itemize}
\item[a)] when there is an exact matching between the string memorized in 
          the searching algorithm (and thus, a string necessarily less 
          than $\log_2 N + k_2$ bits long) and one of the strings provided 
          by the first algorithm.
\item[b)] when there is an exact matching between a string longer than
          $\log_2 N + k_2$ bits, but generated by a sub-program of the 
          searching algorithm (and thus, again less complex than 
          $\log_2 N + k_2$ bits), and one of the strings provided by the 
          first algorithm.
\end{itemize}
In both cases, the searching algorithm is able to find only strings less
complex than $\log_2 N + k_2$ bits.

But we already know that what we are searching is not compressible in less
than $N$ bits. The algorithmic complexity of the desired theorem (or,
equivalently, the algorithmic complexity of its proof, since the proof
provides the theorem) is necessarily greater than $N$ bits because it
generates a string $s$ more complex than $N$ bits as output. Therefore, we
know from the beginning that the searching algorithm will not be
able to spot it. An easy way to show that fact ironically follows from the
very notion of algorithmic complexity, as the Chaitin's theorem does.

Let us forget for a while the subject of this paper. According to the
notion of algorithmic complexity, {\em any} string $l$ of algorithmic
complexity greater than $N$ bits could not be unequivocally recognized
using an algorithm of a complexity of $\log_2 N + k_2$ bits, if $\log_2 N
+ k_2$ is less that $N$. As if it were the case, it would be possible to
write an algorithmic procedure of complexity $\sim \log_2 N + k_2$ which
systematically generates all the strings in increasing size order until it
recognizes the desired one. But, in this case there would be a
contradiction since the algorithmic complexity of the generated string
would necessarily be not greater than $\sim \log_2 N + k_2$ bits. Now,
returning to our point, what difference does it make if $l$ is a whatever
string more complex than $N$ bits or is the string coding our desired
theorem? And this limit affects the {\em searching algorithm}, not the
{\em proof-checking one}. Indeed, all the above does not eliminate {\em a
priori} the possibility that among all the valid proofs generated and
checked by the first algorithm there is one which states that a binary
string $s$ is more complex than $N$ bits.

It is known from the work of A.~M.~Turing that algorithms coding formal
axiomatic systems have inescapable limits in proof-checking. These limits
are strictly related to the undecidability of the {\em halting problem},
namely the non-existence of a general algorithmic procedure to establish
in a finite number of steps whether a computer program halts or not,
simply by analyzing its code. Indeed, it can be demonstrated that the
unsolvability of the halting problem implies the uncomputability of
the program-size complexity.

Instead, in the case of Chaitin~\cite{co} the claimed limit seems to be
not a fundamental one. What can be rigorously deduced from his arguments
is only a trivial limit in the mechanical recognition of the proof, not in
the capability of the system in checking the formal correctness of such
proof.

Substantially, when Chaitin~\cite{co} states that it is possible to write
a program only $\log_2 N +k$ bits long to generate a number of a
complexity greater than $N$ bits, he seems to assume from the beginning 
that it is always possible to compress any theorem of this kind to $k$ bits 
plus the bits which give the binary expansion of $N$. Hence, no wonder if he
gets a logical paradox then.

\section{Conclusions}

From what is outlined above it seems that an obvious, clear relation
between the size of a formal axiomatic system and its
capabilities/limitations in checking proofs can not be claimed. For
instance, a formal axiomatic system of $N$ bits could be able in principle
to check the formal correctness of the proof that a specific sequence is
of a complexity of $C$ bits, with $C\gg N$. The point is that this proof
must be provided as input to the program which codes the proof-checking
system. All this is by no means font of logical paradox. For instance, in
the case such proof exists, and in the case such existence means that the
formal axiomatic system of $N$ bits is able to partially resolve the
halting problem for every string/program in increasing size order up to
$C$ bits, then we have to tell anyway to the formal axiomatic system at
which point to stop its analysis of the halting problem. Namely, we have
to furnish the ordinal number of the last string/program to be checked,
i.e.~an extra $C$ bits\footnote{Let us list all strings of $1,2,3,...,n$
bits, in increasing (lexicographical) order: \begin{center}
\begin{tabular}{lrrrrrrrrrr} $N_\#$ & 1 & 2 & 3 & 4 & 5 & 6 & 7 & 8 & 9 &
...\\ $p$ & 0 & 1 & 00 & 10 & 000 & 100 & 101 & 111 & 0000 & ...  
\end{tabular} \end{center} where $N_\#$ is the ordinal number of the
string, while $p$ denotes the binary string itself. As it can be seen
$p\simeq (N_\#)_2$ and thus the size of $p$ in bits is approximately equal
to $\log_2 N_\#$.}. Otherwise, the system will never halt, vainly trying
to establish the halting problem on the first undecidable program.
Therefore, we simply have a program of $N+C$ bits which generates a string
as complex as $C$ bits.

\section*{Acknowledgments}

I wish to thank Anatoly Vorobey for criticism and useful discussions about 
the subject of this paper.

\end{document}